\documentclass{amsart}

\usepackage{amsmath} 
\usepackage{amssymb}

\newtheorem{theorem}{Theorem}[section] 
\newtheorem{claim}{Claim}[theorem]
\newtheorem{lemma}[theorem]{Lemma} 
\newtheorem{proposition}[theorem]{Proposition} 
\newtheorem{corollary}[theorem]{Corollary} 

\theoremstyle{definition}
\newtheorem{definition}[theorem]{Definition}
\newtheorem{problem}[theorem]{Problem} 

\theoremstyle{remark}

\newtheorem{notation}[theorem]{Notation}
\newtheorem{conclusion}[theorem]{Conclusion}

\numberwithin{equation}{section}
\newcommand{\forces}{\Vdash} 
\newcommand{\bV}{{\bf V}} 
  
\newcommand{\con}{{\mathfrak c}} 
\newcommand{\cf}{{\rm cf}\/} 

\newcommand{\baire}{\omega^{\textstyle \omega}} 
\newcommand{\lh}{\ell g\/} 
\newcommand{\rest}{{\restriction}}
\newcommand{\mrot}{{\rm root}\/} 
\newcommand{\suc}{{\rm succ}} 
\newcommand{\dom}{{\rm dom}} 

\newcommand{\nor}{{\rm {\bf nor}}\/} 
\newcommand{\pos}{{\rm pos}}
\newcommand{\dn}{{\rm dn}}

\newcommand{\basis}{{\bf basis}}
\newcommand{\FC}{{\rm FC}}
\newcommand{\dcl}{{\rm dcl}}
\newcommand{\tree}{{\rm tree}}
\newcommand{\MA}{{\bf MA}}
\newcommand{\FA}{{\bf FA}}

\newcommand{\bC}{{\mathbb C}}
\newcommand{\bD}{{\mathbb D}}
\newcommand{\cF}{{\mathcal F}}
\newcommand{\cG}{{\mathcal G}}
\newcommand{\bH}{{\bf H}}
\newcommand{\cK}{{\mathcal K}}
\newcommand{\bP}{{\mathbb P}}
\newcommand{\cP}{{\mathcal P}}
\newcommand{\bQ}{{\mathbb Q}}
\newcommand{\cS}{{\mathcal S}}
\newcommand{\cU}{{\mathcal U}}

\newcount\skewfactor
\def\mathunderaccent#1#2 {\let\theaccent#1\skewfactor#2
\mathpalette\putaccentunder}
\def\putaccentunder#1#2{\oalign{$#1#2$\crcr\hidewidth
\vbox to.2ex{\hbox{$#1\skew\skewfactor\theaccent{}$}\vss}\hidewidth}}

\newcommand{\gS}{\operatorname{{\mathbf {S}}}_{g, g^\star}}

\renewcommand{\>}{\rangle}

\begin{document}

\title{Towards Martin's Minimum}

\author[T.~Bartoszynski]{Tomek Bartoszynski}
\address{Department of Mathematics and Computer Science\\
 Boise State University\\
 Boise ID 83725, USA} 
\email{tomek@math.idbsu.edu}
\urladdr{http://math.idbsu.edu/$\sim$tomek}
\author[A.~Roslanowski]{Andrzej Ros{\l}anowski}
\address{Department of Mathematics and Computer Science\\
 Boise State University\\
 Boise ID 83725, USA\\
 and Mathematical Institute of Wroclaw University\\
 50384 Wroclaw, Poland} 
\email{roslanow@math.idbsu.edu}
\urladdr{http://math.idbsu.edu/$\sim$roslanow}
\thanks{The second author thanks the KBN (Polish Committee of Scientific
  Research) for partial support through grant 2P03A03114.}


\begin{abstract}
We show that it is consistent with \MA +$\neg$CH that the Forcing
Axiom fails for all forcing notions in the class of $\baire$--bounding
forcing notions with norms of \cite{RoSh:470}.
\end{abstract}

\maketitle

\section{Introduction}
A Forcing Axiom $\FA_\kappa(\bP)$ for a forcing notion $\bP$ is a statement
guaranteeing existence of directed subsets of $\bP$ that meet any member of
a pregiven family of size $\kappa$ of dense subsets of $\bP$. As an axiom,
$\FA_\kappa(\bP)$ is a powerful combinatorial tool that allows one to get
some of the properties of forcing extensions. Even more interesting are
axioms demanding $\FA_\kappa(\bP)$ for all forcing notions in a fixed family
$\cK$. Here the most popular are Martin's Axiom {\MA} postulating
$\FA_\kappa(\bP)$ for each ccc partial order $\bP$ and any cardinal
$\kappa<\con$ and (capturing more forcing notions) Proper Forcing Axiom {\bf
PFA} postulating $\FA_{\aleph_1}(\bP)$ for all proper forcing notions. The
quest for giving the largest possible class of forcing notions $\bP$ for
which the axiom $\FA_{\aleph_1}(\bP)$ may (simultaneously) hold was
successfully accomplished by Foreman, Magidor and Shelah \cite{FMSh:240},
\cite{FMSh:252}, who introduced Martin's Maximum.

In the present paper we want to start investigations in the opposite
directions, looking in some sense for the {\em minimal} version of the
standard Martin's Axiom. That is, we would like to have a model in which
$\neg {\bf CH} + \MA$ holds but $\FA_{\aleph_1}(\bP)$ fails for as many
(necessarily not-ccc) forcing notions as possible. Our attention
concentrates on forcing notions built according to the scheme of {\em norms
on possibilities} of Ros{\l}anowski and Shelah \cite{RoSh:470}. 
(Unfortunately, some familiarity with that paper has to be assumed. In
particular, for all definitions related to norms on possibilities we have to
refer the reader to \cite{RoSh:470}.) 

These lines of investigations have some history already. Steprans proved
that $\neg {\bf CH} + \MA$ is consistent with the negation of the forcing
axiom for the Silver forcing notion (see \cite[\S 2]{CRSW93}). Next, Judah,
Miller and Shelah \cite{JMSh:372}, and Velickovic \cite{Vel91Pos},
 showed that $\neg {\bf CH} + \MA$ does not
imply the forcing axiom for the Sacks forcing notion. It has been a common
believe that the arguments of \cite{JMSh:372} can be repeated for a number
of forcing notions in which conditions are finitely branching trees. For
example, Brendle wrote in the proof of \cite[Proposition 5.1(c)]{Br96a}
(about the method of \cite{JMSh:372}): ``it is easy to see that this
argument works for any forcing notion with compact trees which doesn't have
splitting going on every level''. It seems that \cite{RoSh:470} has provided
the right formalism for specifying which forcing notions can be taken care
of in this context. However, one would like to cover all of
$\baire$--bounding forcing notions from \cite{RoSh:470} (so avoid the
limitations on splittings which in the language of \cite{RoSh:470} would
restrict us to t--omittory trees) and get the failure of $\FA_{\aleph_1}$
for all these forcing notions simultaneously. As remarked by Brendle, the
obvious modifications of the method of \cite{JMSh:372} seem to be not
applicable here. Therefore, we rather follow the Steprans way slightly
generalizing it to be able to deal with a number of forcing notions in the
same model.  
\medskip

\noindent{\bf Notation}\quad Our notation is rather standard and compatible
with that of classical textbooks on Set Theory (like Bartoszynski and
Judah \cite{BaJu95}). However in forcing we keep the older convention that
{\em a stronger condition is the larger one}. 

\begin{notation}
\label{notacja}
\begin{enumerate}
\item For two sequences $\eta,\nu$ we write $\nu\vartriangleleft\eta$ whenever
$\nu$ is a proper initial segment of $\eta$, and $\nu\trianglelefteq\eta$ when
either $\nu\vartriangleleft\eta$ or $\nu=\eta$. The length of a sequence
$\eta$ is denoted by $\lh(\eta)$.
\item The cardinality of the continuum is denoted by $\con$. 
\item For a forcing notion $\bP$, $\Gamma_\bP$ stands for the canonical
$\bP$--name for the generic filter in $\bP$. With this one exception, all
$\bP$--names for objects in the extension via $\bP$ will be denoted with a dot
above (e.g.~$\dot{s}$, $\dot{f}$).
\item Ordinal numbers are denoted by $\alpha,\beta,\gamma,\delta,
\varepsilon, \xi,\zeta$ (with possible indexes); cardinals will be called
$\kappa,\mu$. The first infinite ordinal is $\omega$. The letters $u,v,\eta,
\nu,\rho$ will stand for finite sequences.
\end{enumerate}
\end{notation}

\noindent{\bf Where are the respective definitions in \cite{RoSh:470}?}\quad
As we stated before, we have to assume that the reader is familiar with
\cite{RoSh:470}. (Otherwise we would have to give the list of all needed
definitions and it could be longer than the rest of the paper.) However, for
reader's convenience we list below exact pointers to the descriptions of
the cases of norms on possibilities that are used here.
\begin{itemize}
\item Weak creatures and weak creating pairs: \cite[\S 1.1]{RoSh:470}, in
particular \cite[Definitions 1.1.1, 1.1.3]{RoSh:470}.
\item Creatures and creating pairs: \cite[\S 1.2]{RoSh:470}, in particular
\cite[Definitions 1.2.1, 1.2.2, 1.2.4, 1.2.5]{RoSh:470} (as in
\cite{RoSh:470}, we assume here that all creating pairs are nice and
smooth).
\item Forcing notions $\bQ^*_f(K,\Sigma)$, $\bQ^*_{{\rm w}\infty}$:
\cite[Definitions 1.1.6, 1.1.7, 1.1.10]{RoSh:470}.
\item Tree--creatures and tree--creating pairs: \cite[\S 1.3]{RoSh:470}, in
particular \cite[Definitions 1.3.1, 1.3.3]{RoSh:470}.
\item Forcing notions $\bQ^\tree_0(K,\Sigma)$, $\bQ^\tree_1(K,\Sigma)$:
\cite[Def.\ 1.3.5]{RoSh:470}. 
\item Creating pairs which capture singletons: \cite[Def.\ 2.1.10,
1.2.5(3)]{RoSh:470}.
\item $\bar{2}$--big (tree) creating pairs: \cite[Def.\ 2.2.1,
2.3.2]{RoSh:470}.
\item Halving Property: \cite[Def.\ 2.2.7]{RoSh:470}.
\item $\bH$--fast function $f$: \cite[Def.\ 1.1.12]{RoSh:470}.
\item Simple creating pairs, gluing creating pairs: \cite[Def.\
2.1.7]{RoSh:470}.
\item t--omittory tree--creating pairs: \cite[Def.\ 2.3.4]{RoSh:470}.
\item Strongly finitary creating pairs / tree--creating pairs: \cite[Def.\
1.1.3 + 3.3.4]{RoSh:470}.  
\end{itemize}

\noindent{\bf The forcing notions we want to take care of.}\quad Here we
specify the family of forcing notions for which we want to get the failure
of $\FA_{\aleph_1}$ (with keeping $\MA$).  

\begin{definition}
Let $(K,\Sigma)$ be a weak creating pair for $\bH$. We say that it is {\em
typical} if for each $t\in K$ such that $(\exists u\in\basis(t))(|\pos(u,t)|
=1)$ we have $\nor[t]\leq 1$.
\end{definition}

The reason for the above definition is the following.

\begin{proposition}
\label{divide}
Suppose that $(K,\Sigma)$ is a 2-big creating pair (tree--creating pair,
respectively), $t\in K$ is such that $\nor[t]\geq 4$. Let $u\in\basis(t)$.
Then there are creatures (tree--creatures, resp.) $s_0,s_1\in\Sigma(t)$ such
that $\nor[s_0],\nor[s_1]\geq \nor[t]-2$ and $\pos(u,s_0)\cap\pos(u,s_1)=
\emptyset$. 
\end{proposition}

\begin{proof}
Choose a set $A\subseteq \pos(u,t)$ such that 
\begin{itemize}
\item there is $s^*\in\Sigma(t)$ with $\pos(u,s^*)\subseteq A$ and
$\nor[s^*]\geq\nor[t]-1$, but 
\item for each $a\in A$ and $s\in\Sigma(t)$, if $\pos(u,s)\subseteq A
\setminus\{a\}$ then $\nor[s]<\nor[t]-1$.
\end{itemize}
Clearly it is possible; necessarily $|A|\geq 2$ (remember that $(K,\Sigma)$
is typical). Fix $a\in A$. Applying bigness to $s^*$ we get $s_0\in
\Sigma(s^*)\subseteq\Sigma(t)$ such that $\nor[s_0]\geq\nor[t]-2$ and
$\pos(u,s_0)\subseteq A\setminus\{a\}$. On the other hand, by the choice of
the set $A$ (and bigness) we find $s_1\in\Sigma(t)$ such that $\nor[s_1]\geq
\nor[t]-1$ and $\pos(u,s_1)\subseteq \pos(u,t)\setminus (A\setminus\{a\})$.
\end{proof}

\begin{definition}
\label{classK}
Let $\cK$ be the family of all non-trivial forcing notions of one of the
following types: 
\begin{enumerate}
\item $\bQ^\tree_1(K,\Sigma)$, where $(K,\Sigma)\in {\mathcal H}(\aleph_1)$
is a strongly finitary 2-big typical tree--creating pair for a function
$\bH\in{\mathcal H}(\aleph_1)$;
\item $\bQ^\tree_1(K,\Sigma)$, where $(K,\Sigma)\in {\mathcal H}(\aleph_1)$
is a strongly finitary typical t-omittory tree--creating pair for a function
$\bH\in{\mathcal H}(\aleph_1)$ (note that in this case $\bQ^\tree_1(K,
\Sigma)$ is a dense subforcing of $\bQ^\tree_0(K,\Sigma)$, see
\cite[2.3.5]{RoSh:470}); 
\item $\bQ^*_f(K,\Sigma)$, where $(K,\Sigma)\in {\mathcal H}(\aleph_1)$ is a
strongly finitary typical creating pair for $\bH\in {\mathcal H}(\aleph_1)$,
$f:\omega\times\omega\longrightarrow\omega$ is an $\bH$-fast function, $(K,
\Sigma)$ is $\bar{2}$-big, has the Halving Property and is either simple or
gluing; 
\item $\bQ^*_{{\rm w}\infty}(K,\Sigma)$, where $(K,\Sigma)\in {\mathcal
H}(\aleph_1)$ is a strongly finitary typical creating pair for
$\bH\in{\mathcal H}(\aleph_1)$, $(K,\Sigma)$ captures singletons. 
\end{enumerate}
\end{definition}

\begin{theorem}
[See {\cite[\S 2.3]{RoSh:470}}]
All forcing notions in the class $\cK$ are proper and $\baire$--bounding.
\end{theorem}


Now we may state the main result of this paper.

\begin{theorem}
\label{main}
Assume $\kappa=\kappa^{<\kappa}$. Then there is a ccc forcing notion $\bP$
of size $\kappa$ such that in $\bV^{\bP}$:
\begin{itemize}
\item $\con=\kappa$ and $\MA$ holds true, but
\item $\FA_{\aleph_1}(\bQ)$ fails for every forcing notion $\bQ\in\cK$.
\end{itemize}
\end{theorem}

The class $\cK$ includes all $\baire$--bounding forcing
notions presented in \cite{RoSh:470} (modulo the demand that they are
additionally required to be in ${\mathcal H}(\aleph_1)$). In particular, the
following forcing notions are in $\cK$:
\begin{enumerate}
\item Sacks-like forcings. Let $\{I_n: n \in \omega\}$ be a sequence
  of finite sets. Every perfect set $P \subseteq \prod_n I_n$
  corresponds to a tree $T$. Consider the forcing notion ${\mathbf P}$ 
  which consists of trees $T$ such that 
$$ \forall s \in T \ \exists t \supseteq s \ \suc(t)=I_{|t|},$$ ordered by
inclusion. 
If $I_n=2$ for all $n$ we get the Sacks forcing $\mathbf S$. If
$I_n=k$ ($n \in \omega $) then we get forcing notions $\bD_k$ (for
$k<\omega$) of Newelski and 
Ros{\l}anowski \cite{NeRo93}. Finally, if $I_n=h(n)$ then we get the
Shelah forcing ${\mathbf Q}_h$ of \cite{GJS92}.
\item Silver-like  forcings.
 For $h \in \omega^\omega$ let ${\mathbf P}_h=\{f: \omega \setminus \dom(f)\in
  [\omega]^\omega, \ \forall n \ (n \in \dom(f) \rightarrow f(n)\leq
  h(n))\}$. For $f,g \in
  {\mathbf P}_h$, $f \geq g$ if $g \subseteq f$.
If $h(n)=2$ ($ n \in \omega $) then ${\mathbf P}_h$ is Silver forcing
and for $h(n)=n$ we get Miller forcing from \cite{Mil81Som}.
\item forcing notion $\bQ_{f,g}$ from \cite{BJSh:368} its siblings from 
  \cite{CiSh:653} and \cite{GoShMan93}. Suppose that $f \in
  \omega^\omega$ and $g \in \omega^{\omega \times 
\omega}$ are two functions such that 
\begin{enumerate}
\item $f(n) > \prod_{j<n} f(j)$ for $ n \in \omega$,
\item $g(n,j+1) > f(n) \cdot g(n,j)$ for $n,j \in \omega$, and 
\item $\min\left\{j \in \omega: g(n,j) > f(n+1)\right\}
\stackrel{n \rightarrow \infty}{\longrightarrow} \infty$. 
\end{enumerate}
Define
$Seq^{f} = \{s \in \omega^{<\omega}: \forall j < |s|\ s(j)\leq f(j)\}$
and let
$T \in {\mathbf {PT}}_{f,g}$ if 
\begin{enumerate}
\item $T$ is a perfect subtree of $Seq^{f}$, and 
\item there exists a function $r \in \omega^{\omega}$, $\lim_{n
  \rightarrow \infty} r(n)=\infty$ such that 
$$\forall s \in T \ (stem(T) \subseteq s \rightarrow 
|\suc_{T}(s)| \geq g(|s|,r(|s|))).$$
\end{enumerate}
\item Forcing $\gS$ from \cite{BaJu95}.
Let $g, g^\star \in \omega^\omega$ be  two strictly increasing functions
such that $g(0)=0$, $g^\star(0)=1$ and 
$g(n) << g^\star(n) << g(n+1)$ for all $n>0$.

For $n \in
\omega$ let $$P_n = \left\{a \subseteq g(n+1):
|a|=\frac{g(n+1)}{2^n}\right\}.$$  

For a set $A \subseteq P_n$ define
$$\nor(A)=\min\{|X|: \forall a \in A \ X \not \subseteq a\}.$$
Let 
$$\gS=\left\{\<t_n: n \in \omega\>: \forall n \ t_n \in P_n\ \&\
\forall k\ 
\limsup_{n \rightarrow \infty} \frac{\nor(t_n)}{g^\star(n+1)^k}=\infty\right\}.$$
For $p=\<t^0_n: n \in \omega\>$ and $q=\<t^1_n: n \in \omega\>$ we
define $p \geq q$ if $t^0_n \subseteq t^1_n$ for all $n$.
\item More complicated forcing notions from \cite{FrSh:406} and
  \cite{BaSh607}.  

\end{enumerate}

\section{$\cS$--families of good graphs}
In this section we introduce a property of families of graphs that will be
one of our main tools later.
\begin{definition}
\label{goodgraph}
Let $\cU$ be a countable basis of a topology on a set $X$, $A\subseteq
X\times\omega_1$. A triple $\cG=(A,\cU,E)$ is a {\em good graph} if the
following conditions are satisfied: 
\begin{enumerate}
\item[(a)] $E\subseteq [A]^{\textstyle 2}$,  and
\item[(b)] if $(x_0,\alpha_0), (x_1,\alpha_1)\in A$ are distinct then
$x_0\neq x_1$, and
\item[(c)] if $(x_0,\alpha_0), (x_1,\alpha_1)\in A$, $x_0\neq x_1$ and  
$\{(x_0,\alpha_0),(x_1,\alpha_1)\}\notin E$ then there are disjoint
$U_0,U_1\in\cU$ such that $x_0\in U_0$, $x_1\in U_1$ and 
\[(\forall (x_0',\alpha_0'), (x_1',\alpha_1')\in A)(x_0'\in U_0\ \&\ x_1'\in
U_1\quad\Rightarrow\quad\{(x_0',\alpha_0'), (x_1',\alpha_1')\}\notin E).\]
\end{enumerate}
\end{definition}

\begin{definition}
\label{Sfamily} 
Suppose that $\cF=\{\cG_\zeta:\zeta<\xi\}$ is a family of good graphs,
$\cG_\zeta=(A_\zeta,U_\zeta,E_\zeta)$. 
\begin{enumerate}
\item Let $0<m<\omega$. An $m$--selector for $\cF$ is a set $S\subseteq
(\bigcup\limits_{\zeta<\xi}A_\zeta)^{\textstyle m}$ such that for some 
(not necessarily distinct) $\zeta_0,\ldots,\zeta_{m-1}<\xi$ we have 
\begin{enumerate}
\item[$(\alpha)$] if $\nu\in S$, $\ell<m$ then $\nu(\ell)\in A_{\zeta_\ell}$,
\item[$(\beta)$]  if $\nu,\rho\in S$ are distinct, then for some $\ell<m$ we
have $\{\nu(\ell),\rho(\ell)\}\in E_{\zeta_\ell}$. 
\end{enumerate}
\item $\cF$ is called an {\em $\cS$--$m$-family} if there is no uncountable 
$m$--selector for $\cF$.
\item $\cF$ is an {\em $\cS$--family} (of good graphs) if it is an
$\cS$--$m$-family for each $m<\omega$.
\end{enumerate}
\end{definition}

Let us show how we are going to use $\cS$--families of good graphs.

\begin{definition}
\label{represent}
Let $\bP=(\bP,\leq)$ be a forcing notion and let $\cG=(A,\cU,E)$ be a good
graph (with $A\subseteq X\times\omega_1$, $\cU$ a countable basis of a
topology on $X$). We say that $\cG$ is {\em densely representable by
$\bot_{\bP}$} if there is a one-to-one mapping $\pi:A\longrightarrow \bP$
such that  
\begin{enumerate}
\item[(a)] for each $\alpha<\omega_1$, the set $\{\pi(x,\alpha):
(x,\alpha)\in A\}$ is dense in $\bP$, and
\item[(b)] if $(x_0,\alpha_0),(x_1,\alpha_1)\in A$ are distinct, $\{(x_0,
\alpha_0),(x_1,\alpha_1)\}\notin E$ then the conditions $\pi(x_0,\alpha_0)$,
$\pi(x_1,\alpha_1)$ are incompatible (in $\bP$).
\end{enumerate}
\end{definition}

\begin{proposition}
\label{notFAP}
Let  $\cF$ be an $\cS$--$1$--family of good graphs. Suppose that $\bP$ is a
forcing notion such that some $\cG\in\cF$ is densely representable by
$\bot_{\bP}$. Then $\FA_{\aleph_1}(\bP)$ fails.  
\end{proposition}

\begin{proof}
Let $\cG=(A,\cU,E)$ and let $\pi:A\longrightarrow\bP$ witness that $\cG$ is
densely representable by $\bot_{\bP}$. Let $H_\xi=\{\pi(x,\xi):(x, \xi)\in
A\}$ (for $\xi<\omega_1$). Then the sets $H_\xi$ are dense in $\bP$ by
\ref{represent}(a). We claim that they witness the failure of
$\FA_{\aleph_1}(\bP)$. So suppose that $G\subseteq\bP$ is a directed set
which meets each $H_\xi$. For every $\xi<\omega_1$ choose $(x_\xi,\xi)\in A$
such that $\pi(x_\xi,\xi)\in G\cap H_\xi$. Look at the set $S=\{\langle(
x_\xi,\xi)\rangle:\xi<\omega_1\}$ --- it is an uncountable 1--selector from
$\cF$, a contradiction. 
\end{proof}

One of problems that we will have to take care of when building the forcing
notion needed for \ref{main} is preserving ``being an $\cS$--family''. A
part of this difficulty will be dealt with by ``killing the ccc of bad
forcing notions''.   

\begin{proposition}
\label{notccc}
Let $\cF$ be an $\cS$--family of good graphs. Suppose that $\bP$ is a ccc
forcing notion such that 
\[\forces_{\bP}\mbox{`` $\cF$ is not an $\cS$--family''.}\]
Then there is a ccc forcing notion $\bP^{\cF}$ (of size $\aleph_1$) such
that \[\forces_{\bP^{\cF}}\mbox{`` $\bP$ is not ccc and $\cF$ is an
$\cS$--family''.}\]  
\end{proposition}

\begin{proof}
Let $\cF=\{\cG_\zeta:\zeta<\xi\}$, $\cG_\zeta=(A_\zeta,\cU_\zeta,E_\zeta)$,
$\cU_\zeta$ a basis of a topology on $X_\zeta$, $A_\zeta\subseteq X_\zeta
\times\omega_1$. Assume that some condition in $\bP$ forces that $\cF$ is
not an $\cS$--family. Then we find $m<\omega$, $\zeta_0,\ldots,\zeta_{m-1}<
\xi$, $p\in\bP$ and $\bP$--names $\dot{\nu}_\alpha$ for elements of
$(\bigcup\limits_{\ell<m}A_{\zeta_\ell})^{\textstyle m}$ (for $\alpha<
\omega_1$) such that  
\[\begin{array}{ll}
p\forces_{\bP}&\mbox{`` } (\forall\alpha<\beta<\omega_1)(\dot{\nu}_\alpha
\neq \dot{\nu}_\beta)\quad\&\quad (\forall\alpha<\omega_1)(\forall\ell<m)
(\dot{\nu}_\alpha(\ell)\in A_{\zeta_\ell})\quad\&\\
\ &\ \ (\forall\alpha<\beta<\omega_1)(\exists\ell<m)(\{\dot{\nu}_\alpha(
\ell),\dot{\nu}_\beta(\ell)\}\in E_{\zeta_\ell})\mbox{ ''.}
  \end{array}\]
For each $\alpha<\omega_1$ choose a sequence $\nu_\alpha\in
(\bigcup\limits_{\ell<m}A_{\zeta_\ell})^{\textstyle m}$ and a condition
$p_\alpha\geq p$ such that $p_\alpha\forces\mbox{``} \dot{\nu}_\alpha=
\nu_\alpha$''. (So necessarily $(\forall\alpha<\omega_1)(\forall\ell<m)
(\nu(\ell)\in A_{\zeta_\ell})$.) Let $\bQ$ be the following forcing notion:  
\smallskip

{\bf a condition} in $\bQ$ is a finite set $q\subseteq\omega_1$ such that  
\[(\forall\{\alpha,\beta\}\in [q]^{\textstyle 2})(\forall
\ell<m)(\{\nu_\alpha(\ell),\nu_\beta(\ell)\}\notin E_{\zeta_\ell});\]

{\bf the order} is the inclusion (i.e., $q_1\leq q_2$ if and only if
$q_1\subseteq q_2$). 
\smallskip

\noindent Note that if $\alpha,\beta\in q\in\bQ$, $\alpha\neq\beta$ then the
conditions $p_\alpha,p_\beta$ are incompatible (we will use it in
\ref{Pnotccc}). 
\begin{claim}
\label{trick}
Assume $\{q_\varepsilon:\varepsilon<\omega_1\}\subseteq\bQ$, $q_\varepsilon
=\{\alpha^\varepsilon_i:i<k\}$ (the increasing enumeration; $k<\omega$). 
Then there is $Y\in [\omega_1]^{\textstyle \aleph_1}$ such that
\begin{enumerate}
\item[$(\otimes)_Y$] for each $\varepsilon,\delta\in Y$, if $q_\varepsilon,
q_\delta$ are incompatible in $\bQ$ then there are $i<k$ and $\ell<m$ such
that $\{\nu_{\alpha^\delta_i}(\ell), \nu_{\alpha^\varepsilon_i}(\ell)\}\in
E_{\zeta_\ell}$.
\end{enumerate}
\end{claim}
 
\begin{proof}[Proof of the claim] 
For each $\varepsilon<\omega_1$, $\ell<m$ and distinct $i,j<k$ choose
$U^{\varepsilon,i,j}_\ell\in \cU_{\zeta_\ell}\cup\{*\}$ such that if
$\nu_{\alpha^\varepsilon_i}(\ell)=\nu_{\alpha^\varepsilon_j}(\ell)$ then
$U^{\varepsilon,i,j}_\ell= *$, otherwise $U^{\varepsilon,i,j}_\ell,
U^{\varepsilon,j,i}\in \cU_{\zeta_\ell}$ are such that 
\begin{enumerate}
\item[$(\circledast_1)$] $\nu_{\alpha^\varepsilon_i}(\ell)\in
U^{\varepsilon,i,j}_\ell\times\omega_1$,\quad $\nu_{\alpha^\varepsilon_j}(
\ell)\in U^{\varepsilon,j,i}_\ell\times\omega_1$,\quad $U^{\varepsilon,i,
j}_\ell\cap U^{\varepsilon,j,i}_\ell=\emptyset$, \quad and for all $(x_0,
\alpha_0),(x_1,\alpha_1)\in A_{\zeta_\ell}$, 
\[x_0\in U^{\varepsilon,i,j}_\ell\ \&\ x_1\in U^{\varepsilon,j,i}_\ell\quad
\Rightarrow\quad \{(x_0,\alpha_0),(x_1,\alpha_1)\}\notin E_{\zeta_\ell}.\]
\end{enumerate}
(Why possible? Remember the definition of $\bQ$ and \ref{goodgraph}(c).)
Each $\cU_{\zeta_\ell}$ is countable, so there are $U^{i,j}_\ell\in
\cU_{\zeta_\ell}\cup\{*\}$ and an uncountable set $Y\subseteq\omega_1$ such
that 
\begin{enumerate}
\item[$(\circledast_2)$] $(\forall\varepsilon\in Y)(\forall i,j<k, i\neq j)(
\forall\ell<m)(U^{\varepsilon,i,j}_\ell=U^{i,j}_\ell)$.
\end{enumerate}
Suppose that $\varepsilon,\delta\in Y$ and the conditions $q_\varepsilon,
q_\delta$ are incompatible. It means that there are $i,j<k$ and $\ell<m$
such that $\{\nu_{\alpha^\varepsilon_i}(\ell),\nu_{\alpha^\delta_j}(\ell)\}
\in E_{\zeta_\ell}$. We are going to show that we may demand $i=j$ (what
will finish the proof of the claim). So suppose that $i\neq j$. If
$\nu_{\alpha^\varepsilon_i}(\ell)\neq\nu_{\alpha^\varepsilon_j}(\ell)$, then
(by $(\circledast_2)$) $\nu_{\alpha^\delta_i}(\ell)\neq\nu_{\alpha^\delta_j}
(\ell)$ and (by $(\circledast_1)+(\circledast_2)$) $\nu_{\alpha^\delta_j}(
\ell)\in U^{\delta,i,j}_\ell\times\omega_1=U^{\varepsilon,i,j}_\ell\times
\omega_1$. But applying the last part of $(\circledast_1)$ we may conclude
now that $\{\nu_{\alpha^\varepsilon_i}(\ell),\nu_{\alpha^\delta_j}(\ell)\}
\notin E_{\zeta_\ell}$, a contradiction. So $\nu_{\alpha^\varepsilon_i}(
\ell)=\nu_{\alpha^\varepsilon_j}(\ell)$. But then $\nu_{\alpha^\delta_i}(
\ell)=\nu_{\alpha^\delta_j}(\ell)$ and thus $\{\nu_{\alpha^\varepsilon_i}(
\ell),\nu_{\alpha^\delta_i}(\ell)\}\in E_{\zeta_\ell}$.
\end{proof}

\begin{claim}
\label{Qccc}
$\bQ$ is a ccc forcing notion.
\end{claim}

\begin{proof}[Proof of the claim] Suppose that $\{q_\xi:\xi<\omega_1\}
\subseteq\bQ$ is an antichain in $\bQ$. We may assume that, for some
$k<\omega$, $|q_\xi|=k$ for all $\xi<\omega_1$. Let $q_\xi=\{\alpha^\xi_0,
\alpha^\xi_1,\ldots,\alpha^\xi_{k-1}\}$ be the increasing enumeration. Using
\ref{trick} we may find an uncountable $Y\subseteq\omega_1$ such that for
each distinct $\varepsilon,\delta\in Y$ there are $i_{\varepsilon,\delta}<k$
and $\ell_{\varepsilon,\delta}<m$ with $\{\nu_{\alpha^\delta_{i_{\varepsilon,
\delta}}}(\ell_{i_{\varepsilon,\delta}}),\nu_{\alpha^\varepsilon_{i_{
\varepsilon,\delta}}}(\ell_{i_{\varepsilon,\delta}})\}\in E_{\zeta_\ell}$. 
For each $\varepsilon\in Y$ let $\eta_\varepsilon$ be a sequence of length
$k\cdot m$ such that  
\[\eta_\varepsilon(n)=\nu_{\alpha^\varepsilon_i}(\ell)\quad\mbox{ whenever }
n=i\cdot m+\ell,\ i<k,\ \ell<m.\]
Look at the set $\{\eta_\varepsilon:\varepsilon\in Y\}\subseteq
(\bigcup\limits_{\zeta<\xi} A_\zeta)^{\textstyle k\cdot m}$. It should be
clear that it is an uncountable $k\cdot m$--selector for $\cF$ (the clause
$(\beta)$ of \ref{Sfamily}(1) for $\varepsilon,\delta\in Y$ is witnessed by 
$i_{\varepsilon,\delta}\cdot m+ \ell_{\varepsilon,\delta}$). A contradiction.
\end{proof}

\begin{claim}
\label{Qgood}
$\forces_{\bQ}$ `` $\cF$ is an $\cS$--family of good graphs ''.
\end{claim}

\begin{proof}[Proof of the claim]
The only bad thing that could happen after forcing with $\bQ$ is that an
uncountable $m^*$--selector was added (for some $m^*<\omega$). If so, then
we have $m^*<\omega$, $\bQ$--names $\dot{\eta}_\varepsilon$ (for
$\varepsilon<\omega_1$) and a condition $q\in\bQ$ such that 
\[q\forces_{\bQ}\mbox{`` }\{\dot{\eta}_\varepsilon:\varepsilon<\omega_1\}
\mbox{ is an $m^*$--selector for $\cF$,\ \ $\dot{\eta}_\varepsilon$'s are
pairwise distinct ''.}\]
Clearly we may require that for some $\zeta_0,\ldots,\zeta_{m^*-1}$, for each
$\varepsilon<\omega_1$ the condition $q$ forces that $\dot{\eta}_\varepsilon
(\ell)\in A_{\zeta_\ell}$ (for all $\ell<m^*$). For each $\varepsilon<
\omega_1$ pick a sequence $\eta_\varepsilon$ and a condition $q_\varepsilon
\geq q$ such that $q_\varepsilon\forces_\bQ$``$\dot{\eta}_\varepsilon=
\eta_\varepsilon$''. Next, choose an uncountable set $Y\subseteq\omega_1$
and $k<\omega$ such that for $\varepsilon\in Y$, $q_\varepsilon=
\{\alpha^\epsilon_i:i<k\}$ (the increasing enumeration) and $(\otimes)_Y$ of
\ref{trick} holds (possible by \ref{trick}). Let $\rho_\varepsilon$ (for
$\varepsilon\in Y$) be sequences of length $k\cdot m+m^*$ defined by 
\[\rho_\varepsilon(n)=\left\{
\begin{array}{ll}
\nu_{\alpha^\varepsilon_i}(\ell)&\mbox{if }n=i\cdot m+\ell,\ i<k,\ \ell<m,\\
\eta_\varepsilon(\ell)          &\mbox{if }n=k\cdot m+\ell,\ \ell<m^*.
\end{array}\right.\]
Note that for distinct $\varepsilon,\delta\in Y$ we have:
\begin{enumerate}
\item[$(\boxtimes_1)$] if $q_\varepsilon,q_\delta$ are compatible in $\bQ$,
$\varepsilon<\delta<\omega_1$, then for some $\ell<m^*$ we have
$\{\rho_\varepsilon(km+\ell),\rho_\delta(km+\ell)\}=\{\eta_\varepsilon(\ell),
\eta_\delta(\ell)\}\in E_{\zeta_\ell}$;  
\item[$(\boxtimes)_2$] if $q_\varepsilon,q_\delta$ are incompatible in $\bQ$
then there are $i<k$ and $\ell<m$ so that $\{\rho_\varepsilon(im+\ell),
\rho_\delta(im+\ell)\}=\{\nu_{\alpha^\varepsilon_i}(\ell),\nu_{
\alpha^\delta_i}(\ell)\}\in E_{\zeta_\ell}$. 
\end{enumerate}
(Why? $(\boxtimes_1)$ follows from the choice of $\dot{\eta}_\varepsilon$,
$q_\varepsilon$, $(\boxtimes_2)$ is a consequence of $(\otimes)_Y$.) Hence,
$\{\rho_\varepsilon:\varepsilon\in Y\}$ is an uncountable $km+m^*$--selector
from $\cF$, a contradiction. 
\end{proof}

\begin{claim}
\label{Pnotccc}
For some $q\in\bQ$ we have
\[q\forces_{\bQ}\mbox{ `` $\bP$ does not satisfy the ccc ''.}\]
\end{claim}

\begin{proof}[Proof of the claim] 
As we stated before, if $q\in\bQ$ and $\alpha,\beta\in q$ are distinct then
the conditions $p_\alpha,p_\beta$ are incompatible in $\bP$. Since, by
\ref{Qccc}, the forcing notion $\bQ$ is ccc, there is a condition $q\in\bQ$
such that
\[q\forces_{\bQ}\mbox{`` }\{\alpha<\omega_1: \{\alpha\}\in\Gamma_\bQ\}\mbox{
is uncountable ''}.\]
It should be clear now that the condition $q$ forces (in $\bQ$) that $\bP$
fails the ccc.
\end{proof}

Let $q\in\bQ$ be as guaranteed by \ref{Pnotccc} and let $\bP^\cF$ be the
$\bQ$ restricted to elements stronger than that $q$. It should be clear that
$\bP^\cF$ is as required in the proposition.
\end{proof}

\begin{conclusion}
\label{conclusion}
If $\cF$ is an $\cS$--family of good graphs and $\bP$ is a forcing notion
with the Knaster property, then
\[\forces_{\bP}\mbox{`` $\cF$ is an $\cS$--family ''.}\]
\end{conclusion}

\section{Where are our $\cS$--families from?}
It follows from \ref{notFAP} that to make sure that $\FA_{\aleph_1}(\bP)$
fails for all $\bP\in \cK$ it is enough to have an $\cS$--family of good
graphs such that for every $\bP\in\cK$ some $\cG\in\cF$ is densely
representable by $\bot_{\bP}$. In this section we will (almost) show how the
respective $\cS$--family is created in our model. Basically, it will come
from Cohen reals, but interpreted in a special way.

Let $\bP\in\cK$ and let $p$ be a condition in $\bP$. Considering all
possible cases (of \ref{classK}) we define a countable basis $\cU(\bP)$ of a
topology on $X(\bP)\stackrel{\rm def}{=}\bP$, a set $E(\bP)\subseteq
[X(\bP)\times\omega_1]^{\textstyle 2}$ and a forcing notion
$\bC(\bP,p)$. (In the last case we will assume additionally that the
condition $p$ has some special form, which will restrict us to a dense
subset of $\bP$.)  
\medskip
  
\noindent{\sc Case 1}:\quad $\bP=\bQ^\tree_1(K,\Sigma)$, where
$(K,\Sigma)\in {\mathcal H}(\aleph_1)$ is a strongly finitary 2-big typical
tree--creating pair for a function $\bH\in{\mathcal H}(\aleph_1)$.

\noindent For $q=\langle t^q_\eta:\eta\in T^q\rangle\in\bP$ and $N<\omega$
let $U(q,N)$ be the family of all $q'\in\bP$ such that 
\[\mrot(q')=\mrot(q)\quad\mbox{ and }\quad (\forall\eta\in T^{q'})\big(\lh
(\eta)<N\quad\Rightarrow\quad (\eta\in T^q\ \&\ t^q_\eta=t^{q'}_\eta)
\big).\]
Let $\cU(\bP)$ consist of all nonempty sets $U(q,N)$ (for $q\in \bP$,
$N<\omega$). Clearly $\cU(\bP)$ is a countable basis of a topology on $\bP$
(which is the natural product topology). Next we define
\[\begin{array}{ll}
E(\bP)=\{\{(q_0,\alpha_0),(q_1,\alpha_1)\}\in [\bP\!\times\!\omega_1]^{
\textstyle 2}:& \mbox{for each }N<\omega\\
\ &\dcl(T^{q_0})\cap \prod\limits_{i<N}\bH(i)\cap\dcl(T^{q_1})\neq\emptyset\}
  \end{array}\]
(where $\dcl(T)$ is the downward closure of a quasi tree $T$; see
\cite[Def.\ 1.3.1]{RoSh:470}). The forcing notion $\bC(\bP,p)$ is defined as
follows:
\smallskip

\noindent{\bf a condition $r$ in $\bC(\bP,p)$} is a finite system $r=
\langle s^r_\eta: \eta\in \hat{S^r}\rangle$ such that $S^r\subseteq T^p$ is
a (finite) quasi tree, $\mrot(S^r)=\mrot(T^p)$,
\[(\forall\eta\in \hat{S^r})(s^r_\eta\in\Sigma(t^p_\eta)\ \&\
\nor[s^r_\eta]\geq\nor[t^p_\eta]-2\ \&\ \pos(s^r_\eta)=\suc_{S^r}(\eta)),\]
and if $\nor[t^p_\eta]\leq 4$, $\eta\in \hat{S^r}$ then $s^r_\eta=t^p_\eta$;

\noindent{\bf the order of $\bC(\bP,p)$} is the end extension, i.e., $r_0\leq
r_1$ if and only if $S^{r_0}\subseteq S^{r_1}$ and $(\forall\eta\in
\hat{S^{r_0}})(s^{r_0}_\eta=s^{r_1}_\eta)$.
\smallskip

\noindent It should be clear that $\bC(\bP,p)$ is a countable atomless
(remember \ref{divide}) forcing notion, so it is equivalent to the Cohen
forcing. Moreover, the forcing with $\bC(\bP,p)$ adds a condition
$\dot{p}^*=\bigcup\Gamma_{\bC(\bP,p)}\in\bP$ stronger than $p$.
\medskip

\noindent{\sc Case 2}:\quad $\bP=\bQ^\tree_1(K,\Sigma)$, where $(K,\Sigma)
\in {\mathcal H}(\aleph_1)$ is a strongly finitary typical t-omittory
tree--creating pair for a function $\bH\in{\mathcal H}(\aleph_1)$. 

\noindent $\cU(\bP)$ and $E(\bP)$ are defined like in the Case 1. The forcing
notion $\bC(\bP,p)$ is defined similarly too, though we now make an
advantage from ``t--omittory'':
\smallskip

\noindent{\bf a condition $r$ in $\bC(\bP,p)$} is a finite system $r=
\langle s^r_\eta: \eta\in \hat{S^r}\rangle$ such that $S^r\subseteq T^p$ is
a quasi tree, $\mrot(S^r)=\mrot(T^p)$, and for each $\eta\in \hat{S^r}$ 
\begin{itemize}
\item there is a (finite) quasi tree $T^*_\eta\subseteq T^p$ such that
$\mrot(T^*_\eta)= \eta$ and $s^r_\eta\in\Sigma(t^p_\nu:\nu\in\hat{T^*_\eta
})$,
\item $\nor[s^r_\eta]\geq\min\{\nor[t^p_\nu]-2:\nu\in T^*_\eta\}$, 
\end{itemize}
and if $\nor[t^p_\eta]\leq 4$, $\eta\in \hat{S^r}$ then $s^r_\eta=t^p_\eta$;

\noindent{\bf the order of $\bC(\bP,p)$} is the end extension, so $r_0\leq
r_1$ if and only if $S^{r_0}\subseteq S^{r_1}$ and $(\forall\eta\in
\hat{S^{r_0}})(s^{r_0}_\eta=s^{r_1}_\eta)$.
\smallskip

\noindent Again, $\bC(\bP,p)$ is a countable atomless forcing notion. The
forcing with $\bC(\bP,p)$ adds a condition $\dot{p}^*=\bigcup\Gamma_{\bC(
\bP,p)}\in\bP$ stronger than $p$. 
\medskip

\noindent{\sc Case 3}:\quad $\bP=\bQ^*_f(K,\Sigma)$, where $(K,\Sigma)\in
{\mathcal H}(\aleph_1)$ is a strongly finitary typical creating pair for
$\bH\in {\mathcal H}(\aleph_1)$, $f:\omega\times\omega\longrightarrow\omega$
is an $\bH$-fast function, $(K,\Sigma)$ is $\bar{2}$-big, has the Halving
Property and is either simple or gluing.

\noindent For a condition $q=(w^q,t^q_0,t^q_1,\ldots)\in\bP$ and $N<\omega$
we let $U(q,N)$ be the collection of all conditions $q'\in\bP$ such that
$w^q=w^{q'}$ and $t^q_n=t^{q'}_n$ for all $n<N$. Next, $\cU(\bP)$ is the
family of all non-empty sets $U(q,N)$ (for $q\in\bP$ and $N<\omega$). Like
before, $\cU(\bP)$ is a countable basis of the standard (product) topology
on $\bP$. We define  
\[\begin{array}{ll}
E(\bP)&=\{\{(q_0,\alpha_0),(q_1,\alpha_1)\}\in [\bP\times\omega_1]^{
\textstyle 2}: \mbox{ for each } N<\omega\\
\ &\ \ \dcl(\pos(w^{q_0},t^{q_0}_0,\ldots,t^{q_0}_N))\cap \prod\limits_{i<N}
\bH(i)\cap\dcl(\pos(w^{q_1},t^{q_1}_0,\ldots,t^{q_1}_N))\neq\emptyset\}.
  \end{array}
\]
The forcing notion $\bC(\bP,p)$ is defined so that
\smallskip

\noindent{\bf a condition $r$ in $\bC(\bP,p)$} is a finite sequence $r=(w^r,
s^r_0,\ldots,s^r_{n_r})\in\FC(K,\Sigma)$ such that $w^r=w^p$, $s^r_k\in
\Sigma(t^p_k)$, and 
\begin{itemize}
\item if $\nor[t^p_k]>f(n,m^{t^p_k}_\dn)$, $n\geq 4$ then $\nor[s^r_k]>
f(n-2,m^{s^r_k}_\dn)$, 
\item if $\nor[t^p_k]\leq f(4,m^{t^p_k}_\dn)$ then $s^r_k=t^p_k$;
\end{itemize}

\noindent{\bf the order of $\bC(\bP,p)$} is the extension, i.e., $r_0\leq
r_1$ if and only if $n_{r_0}\leq n_{r_1}$ and $s^{r_0}_k=s^{r_1}_k$ for all
$k\leq n_{r_0}$.
\smallskip

\noindent Clearly $\bC(\bP,p)$ is countable and atomless (remember
\ref{divide}). It adds a condition $\dot{p}^*=\bigcup\Gamma_{\bC(\bP,p)}\in
\bP$ stronger than $p$.  
\medskip

\noindent{\sc Case 4}:\quad $\bP=\bQ^*_{{\rm w}\infty}(K,\Sigma)$,
where $(K,\Sigma)\in {\mathcal H}(\aleph_1)$ is a strongly finitary typical
creating pair for $\bH\in{\mathcal H} (\aleph_1)$, $(K,\Sigma)$ captures
singletons.  

\noindent Both $\cU(\bP)$ and $E(\bP)$ are defined like in Case 3. As we
said before, here we will require that $p$ is of special form. Namely, we
demand that $p=(w^p,t^p_0,t^p_2,\ldots)\in\bP$ is such that for some
strictly increasing sequence $\langle m_k:k<\omega\rangle\subseteq\omega$,
$m_0=0$ and for each $k\in\omega$:
\begin{itemize}
\item $\nor[t^p_{m_k}]\geq 4+k$, and
\item if $m_k+1<m_{k+1}$ then for some (equivalently: all) $u\in\pos(w^p,
t^p_0,\ldots,t^p_{m_k})$ we have $|\pos(u,t^p_{m_k+1},\ldots, t^p_{m_{k+1}
-1})|=1$. 
\end{itemize}
(Since $(K,\Sigma)$ captures singletons the conditions of this form are
dense in $\bP$.) Now, the forcing notion $\bC(\bP,p)$ is defined so that 
\smallskip

\noindent{\bf a condition $r$ in $\bC(\bP,p)$} is a finite sequence $r=(w^r,
s^r_0,\ldots,s^r_{n_r})\in\FC(K,\Sigma)$ such that $w^r=w^p$, $s^r_\ell\in
\Sigma(t^p_\ell)$, and 
\[(\forall k<\omega)(m_{2k}\leq n_r\ \Rightarrow\ \nor[s^r_{m_{2k}}]\geq
\nor[t^p_{m_{2k}}]-2;\] 

\noindent{\bf the order of $\bC(\bP,p)$} is the extension, i.e., $r_0\leq
r_1$ if and only if $n_{r_0}\leq n_{r_1}$ and $s^{r_0}_k=s^{r_1}_k$ for all
$k\leq n_{r_0}$.
\smallskip

\noindent Again, $\bC(\bP,p)$ is countable and atomless, and it adds a
condition $\dot{p}^*=\bigcup\Gamma_{\bC(\bP,p)}\in\bP$ stronger than $p$.  

\begin{lemma}
\label{oneforce}
Suppose $\bP\in\cK$ and $p\in\bP$ is such that $\bC(\bP,p)$ is defined. Let
$r\in\bC(\bP,p)$. Then there are two conditions $r_0,r_1\in\bC(\bP,p)$
stronger than $r$ and basic open sets $U_0,U_1,U\in \cU(\bP)$ such that
\begin{itemize}
\item $U_0\cap U_1=\emptyset$, $p\in U$,
\item if $p'\in U$ and $\bC(\bP,p')$ is defined then $r,r_0,r_1\in\bC(\bP,
p')$, $r_0,r_1$ are stronger than $r$ (in $\bC(\bP,p')$) and
$r_0\forces_{\bC(\bP,p')}\dot{p}^*\in U_0$, $r_1\forces_{\bC(\bP,p')} 
\dot{p}^*\in U_1$, 
\item $(\forall q_0\in U_0)(\forall q_1\in U_1)(\forall\alpha_0,\alpha_1<
\omega_1)(\{(q_0,\alpha_0),(q_1,\alpha_1)\}\notin E(\bP))$. 
\end{itemize}
(Note that these formulas are absolute.)
\end{lemma}

\begin{proof}
In Cases 1 and 3 (of \ref{classK}) use \ref{divide}; in other cases use
directly the assumption that $(K,\Sigma)$ is typical.
\end{proof}

The $\cS$--families in our model will be created by choosing sets
$A\subseteq X(\bP)\times \omega_1$ (for each $\bP\in\cK$) so that in each
pair $(q,\alpha)\in A$ the first coordinate $q$ is added generically by the
forcing notion $\bC(\bP,p)$ (for some condition $p\in\bP$). For this we will
use the finite support product of $\aleph_1$ copies of $\bC(\bP,p)$, denoted
by $\bC^{\omega_1}(\bP,p)=\prod\limits_{\delta<\omega_1}\bC(\bP,p)$ (so
a condition in $\bC^{\omega_1}(\bP,p)$ is a finite function $c:\dom(c)
\longrightarrow\bC(\bP,p)$ and the order is the natural one). The forcing
with $\bC^{\omega_1}(\bP,p)$ adds the set   
\[\dot{Z}^\bP_p=\big\{(q,\alpha):\alpha<\omega_1\ \&\ q=\bigcup\{r:
(\alpha,r)\in\Gamma_{\bC^{\omega_1}(\bP,p)}\}\big\}.\]
(Sets of these form will be used to build a good graph $\cG$ densely
representable by $\bot_{\bP}$.)

\begin{definition}
\label{isomorph}
Suppose that $\bP\in\cK$, $p\in\bP$ and $\bC(\bP,p)$ is defined. We say that
two conditions $\bar{c}_0,\bar{c}_1\in\bC^{\omega_1}(\bP,p)$ are {\em
isomorphic\/} (and then we write $\bar{c}_0\sim \bar{c}_1$) if $|\dom(
\bar{c}_0)|=|\dom(\bar{c}_1)|$ and if $H:\dom(\bar{c}_0)\longrightarrow
\dom(\bar{c}_1)$ is the order preserving bijection then $\bar{c}_1(H(\alpha
))=\bar{c}_0(\alpha)$ for each $\alpha\in\dom(\bar{c}_0)$. 
\end{definition}

\noindent (Note that there are countably many isomorphism types of
conditions in $\bC^{\omega_1}(\bP,p)$.) 

The main technical advantage of using the forcing notions $\bC^{\omega_1}
(\bP,p)$ to create our $\cS$--families is presented by the following lemma.

\begin{lemma}
\label{maintech}
Let $\bP^0,\ldots,\bP^k\in\cK$. Suppose that $\langle \bP_\xi,\dot{\bQ}_\xi:
\xi<\gamma\rangle$ is a finite support iteration of ccc forcing
notions, $\gamma$ is a limit ordinal. Furthermore assume that for some
disjoint sets $I_0,\ldots,I_k\subseteq\gamma$ we have
\begin{enumerate}
\item[$(\alpha)$] if $\xi\in I_\ell$ then $\forces_{\bP_\xi}$``
$\dot{\bQ}_\xi=\bC^{\omega_1}(\bP^\ell,p^\xi)$ for some $\dot{p}^\xi\in
\bP^\ell$'', 
\item[$(\beta)$] for $\zeta\leq\gamma$, $\dot{A}^\ell_\zeta$ is the
$\bP_\zeta$--name for the set $\bigcup\{\dot{Z}^{\bP^\ell}_{\dot{p}^\xi}:\xi
\in I_\ell\cap\zeta\}\subseteq\bP^\ell\times \omega_1$, 
\item[$(\gamma)$] for each $\zeta<\gamma$
\[\forces_{\bP_\zeta}\mbox{`` }\{(\dot{A}^\ell_\zeta, \cU(\bP^\ell),
E(\bP^\ell)\cap [\dot{A}^\ell_\zeta]^{\textstyle 2}): \ell\leq k\}\mbox{ is
an $\cS$--family of good graphs ''.}\] 
\end{enumerate}
Then 
\[\forces_{\bP_\gamma}\mbox{`` }\{(\dot{A}^\ell_\gamma, \cU(\bP^\ell),
E(\bP^\ell)\cap [\dot{A}^\ell_\gamma]^{\textstyle 2}): \ell\leq k\}\mbox{ 
is an $\cS$--family of good graphs. ''}\]
\end{lemma}

\begin{proof}
For a condition $q_0\in\bP_\gamma$ we may find a stronger condition $q\in
\bP_\gamma$ with the following property 
\begin{enumerate}
\item[$(*)_q$] for each $\xi\in I_\ell\cap\dom(q)$, $\ell\leq k$ there
are $\bar{c}(q,\xi)$, $\bar{c}_0(q,\xi)$, $\bar{c}_1(q,\xi)$, and
$U(q,\xi)$, $\bar{U}_0(q,\xi)$, $\bar{U}_1(q,\xi)$ (objects, not names) such
that the condition $q\rest\xi$ forces the following: 
\begin{itemize}
\item $U(q,\xi)\in\cU(\bP^\ell)$, $\dot{p}_\xi\in U(q,\xi)$, $q(\xi)=
\bar{c}(q,\xi)\in\bC^{\omega_1}(\bP^\ell,\dot{p}_\xi)$,
\item $\bar{U}_0(q,\xi),\bar{U}_1(q,\xi):\dom(\bar{c}(q,\xi))\longrightarrow
\cU(\bP^\ell)$, $\bar{U}_0(q,\xi)(\varepsilon)\cap\bar{U}_1(q,\xi)(
\varepsilon)=\emptyset$ for each $\varepsilon\in\dom(\bar{c}(q,\xi))$, 
\item for each $p'\in U(q,\xi)$ such that $\bC(\bP^\ell,p')$ is defined: 
$\bar{c}(q,\xi),\bar{c}_0(q,\xi),\bar{c}_1(q,\xi)$ are in $\bC^{\omega_1}(
\bP^\ell,p')$, the conditions $\bar{c}_0(q,\xi),\bar{c}_1(q,\xi)$ are stronger
than $\bar{c}(q,\xi)$ and $\dom(\bar{c}_0(q,\xi))=\dom(\bar{c}_1(q,\xi))=
\dom(\bar{c}(q,\xi))$, and 

if $i<2$, $\varepsilon\in\dom(\bar{c}(q,\xi))$, and $\dot{p}^*_\varepsilon$
is the $\bC^{\omega_1}(\bP^\ell,p')$--name for the $\varepsilon^{\rm th}$
generic real (i.e., $\bigcup\{r:(\alpha,r)\in\Gamma_{\bC^{\omega_1}(
\bP^\ell,p')}\}$) then 
\[\bar{c}_i(q,\xi)\forces_{\bC^{\omega_1}(\bP,p')}\mbox{`` }
\dot{p}^*_\varepsilon\in \bar{U}_i(q,\xi)(\varepsilon)\mbox{ '',}\] 
\item for each $\varepsilon\in\dom(\bar{c}(q,\xi))$ and $q_0\in\bar{U}_0(q,
\xi)(\varepsilon)$, $q_1\in\bar{U}_1(q,\xi)(\varepsilon)$ we have
\[(\forall\varepsilon_0,\varepsilon_1<\omega_1)(\{(q_0,\varepsilon_0),(q_1,
\varepsilon_1)\}\notin E(\bP^\ell)).\]    
\end{itemize}
\end{enumerate}
[Why? Just apply \ref{oneforce} (and remember that supports are finite).]
>From now on we will restrict ourselves to conditions $q\in\bP_\gamma$ with
the property $(*)_q$ (what is allowed as they are dense in $\bP_\gamma$). So
we will assume that for each condition $q$ under considerations and $\xi\in
I_\ell\cap\dom(q)$, $\ell\leq k$, the objects (not names!) $\bar{c}(q,\xi)$,
$\bar{c}_0(q,\xi)$, $\bar{c}_1(q,\xi)$, $U(q,\xi)$, $\bar{U}_0(q,\xi)$,
$\bar{U}_1(q,\xi)$ are defined and have the respective properties.   

Note that, in $\bV^{\bP_\gamma}$, if $\ell\leq k$, $\xi_0,\xi_\ell\in
I_\ell$, $(q,\alpha_0)\in \dot{Z}^{\bP^\ell}_{\dot{p}^{\xi_0}}$ and $(q,
\alpha_1)\in\dot{Z}^{\bP^\ell}_{\dot{p}^{\xi_1}}$ then $\alpha_0=\alpha_1$
and $\xi_0=\xi_1$ (remember \ref{oneforce}). Therefore, we may label
elements of $\dot{A}^\ell_\gamma$ by pairs from $I_\ell\times\omega_1$ and
allow ourselves small abuse of notation identifying $(\xi,\alpha)\in 
I_\ell\times\omega_1$ with the respective $(q,\alpha)\in
\dot{Z}^{\bP^\ell}_{\dot{p}^\xi}$. Next let $E_\ell=E(\bP^\ell)\cap
[\dot{A}^\ell_\gamma]^{\textstyle 2}$. 

Now, suppose that some condition $q'\in\bP_\gamma$ forces that 
\[\{(\dot{A}^\ell_\gamma, \cU(\bP^\ell),E_\ell): \ell\leq k\}\mbox{ is not an
$\cS$--family.}\]
Then we may find a condition $q\in\bP_\gamma$, an integer $m<\omega$,
$\ell_0,\ldots,\ell_{m-1}\leq k$ (not necessarily distinct) and
$\bP_\gamma$--names $\dot{\nu}_\alpha$ (for $\alpha<\omega_1$) of sequences
of length $m$ such that 
\[\begin{array}{ll}
q\forces_{\bP_\gamma}&\mbox{`` }(\forall\alpha<\beta<\omega_1)(
\dot{\nu}_\alpha\neq\dot{\nu}_\beta)\ \&\ (\forall\alpha<\omega_1)(\forall i<
m)(\dot{\nu}_\alpha(i)\in I_{\ell_i}\times\omega_1)\\
\ &\ \ (\forall\alpha<\beta<\omega_1)(\exists i<m)(\{\dot{\nu}_\alpha(i),
\dot{\nu}_\beta(i)\}\in E_{\ell_i})\mbox{ ''.}
  \end{array}\]
For each $\alpha<\omega_1$ choose a condition $q_\alpha\in\bP_\gamma$
(satisfying $(*)_{q_\alpha}$ and) stronger than $q$ and a sequence
$\nu_\alpha\in\prod\limits_{i<m}(I_{\ell_i}\times\omega_1)$ such that
$q_\alpha\forces\dot{\nu}_\alpha=\nu_\alpha$ and 
\[(\forall i<m)(\nu_\alpha(i)=(\xi,\varepsilon)\ \Rightarrow\
\xi\in\dom(q_\alpha)\ \&\ \varepsilon\in\dom(q_\alpha(\xi))).\]
Now we consider two cases.
\smallskip 

\noindent{\sc Case A:}\qquad $\cf(\gamma)\neq\omega_1$.\\
Then for some $\zeta<\gamma$, for uncountably many $\alpha<\omega_1$,
$\dom(q_\alpha)\subseteq\zeta$. Let $G\subseteq\bP_\gamma$ be a generic over
$\bV$ and work in $\bV[G\cap\bP_\zeta]$. Because of the ccc of $\bP_\zeta$,
the set $\{\alpha<\omega_1: q_\alpha\in G\cap\bP_\zeta\}$ is uncountable, so
we get an uncountable $m$--selector from $\{(\dot{A}^\ell_\zeta,
\cU(\bP^\ell), E(\bP^\ell)\cap [\dot{A}^\ell_\zeta]^{\textstyle 2})^{G\cap
\bP_\zeta}:\ell\leq k\}$ (in $\bV[G\cap\bP_\zeta]$), contradicting the
assumption $(\gamma)$.  
\medskip

\noindent{\sc Case B:}\qquad $\cf(\gamma)=\omega_1$.\\
If for some $\zeta<\gamma$ the set $\{\alpha<\omega_1:\dom(q_\alpha)
\subseteq \zeta\}$ is uncountable then we may repeat the arguments of Case
A. So assume that $\{\alpha<\omega_1:\dom(q_\alpha)\subseteq\zeta\}$ is
countable for each $\zeta<\gamma$.

Applying ``standard cleaning procedure'' and passing to an uncountable
subsequence (and possibly increasing our conditions) we may assume that
$|\dom(q_\alpha)|=N$ for each $\alpha<\omega_1$ and, letting
$\{\xi^\alpha_0,\ldots,\xi^\alpha_{N-1}\}$ be the increasing enumeration of
$\dom(q_\alpha)$:    
\begin{enumerate}
\item $\{\dom(q_\alpha):\alpha<\omega_1\}$ forms a $\Delta$--system with
heart $u^*$,
\item for some $n^*<N$ and $\zeta^*<\gamma$, we have 
$(\forall\alpha<\omega_1)(\forall j<n^*)(\xi^\alpha_j<\zeta^*)$ and 
$(\forall\alpha<\beta<\omega_1)(\zeta^*<\xi^\alpha_{n^*}\leq\xi^\alpha_{N-1}
<\xi^\beta_{n^*})$ (so necessarily $u^*\subseteq\zeta^*$),
\item $\sup\{\xi^\alpha_{n^*}:\alpha<\omega_1\}=\gamma$,
\item $(\forall\alpha,\beta<\omega_1)(\forall\ell\leq k)(\forall j<N)(
\xi^\alpha_j\in I_\ell\ \Leftrightarrow\ \xi^\beta_j\in I_\ell)$,
\item if $\alpha,\beta<\omega_1$, $\ell\leq k$, $j<N$ and $\xi^\alpha_j\in
I_\ell$ then $U(q_\alpha,\xi^\alpha_j)=U(q_\beta,\xi^\beta_j)$,
$\bar{c}(q_\alpha,\xi^\alpha_j)\sim \bar{c}(q_\beta,\xi^\beta_j)$,
$\bar{c}_0(q_\alpha,\xi^\alpha_j)\sim \bar{c}_0(q_\beta,\xi^\beta_j)$,
$\bar{c}_1(q_\alpha,\xi^\alpha_j)\sim \bar{c}_1(q_\beta,\xi^\beta_j)$ (see
\ref{isomorph}), and 
\begin{enumerate}
\item[$(*)$] if $H:\dom(\bar{c}(q_\alpha,\xi^\alpha_j))\longrightarrow\dom
(\bar{c}(q_\beta,\xi^\beta_j))$ is the order preserving bijection then $H$
is the identity on $\dom(\bar{c}(q_\alpha,\xi^\alpha_j))\cap\dom(\bar{c}(
q_\beta,\xi^\beta_j))$ and for each $\varepsilon\in\bar{c}(q_\alpha,
\xi^\alpha_j)$ 
\[\begin{array}{l}
\bar{U}_0(q_\alpha,\xi^\alpha_j)(\varepsilon)=\bar{U}_0(q_\beta,
\xi^\beta_j)(H(\varepsilon)),\quad \bar{U}_1(q_\alpha,\xi^\alpha_j)(
\varepsilon)=\bar{U}_1(q_\beta,\xi^\beta_j)(H(\varepsilon))\quad\mbox{ and}\\
(\forall i<m)(\nu_\alpha(i)=(\xi^\alpha_j,\varepsilon)\ \Leftrightarrow\
\nu_\beta(i)=(\xi^\beta_j,H(\varepsilon))). 
  \end{array}\]  
\end{enumerate}
\end{enumerate}

Let $w^*$ be the set of these $i<m$ that for some (equivalently: all)
$\alpha<\omega_1$ we have $\nu_\alpha(i)\in\zeta^*\times\omega_1$.

\begin{claim}
\label{clx1}
There are $q^*\in\bP_{\zeta^*}$ and $\alpha<\beta<\omega_1$ such that $q^*$
is stronger than both $q_\alpha\rest\zeta^*$ and $q_\beta\rest\zeta^*$ and
$(\forall i\in w^*)(\{\nu_\alpha(i),\nu_\beta(i)\}\notin E_{\ell_i})$.
\end{claim}

\begin{proof}[Proof of the claim] 
Let $G_{\zeta^*}\subseteq\bP_{\zeta^*}$ be a generic filter over $\bV$. Work
in $\bV[G_{\zeta^*}]$. By the ccc of $\bP_{\zeta^*}$, the set $\{\alpha<
\omega_1:q_\alpha\rest\zeta^*\in G_{\zeta^*}\}$ is uncountable. Look at the
sequence $\langle\nu_\alpha\rest w^*:q_\alpha\rest\zeta^*\in G_{\zeta^*}
\rangle$. By assumption $(\gamma)$ of the lemma, it cannot be a
$|w^*|$--selector, so there are $\alpha<\beta<\omega_1$ such that
\[q_\alpha\rest\zeta^*\quad\&\quad q_\beta\rest\zeta^*\quad\&\quad (\forall
i\in w^*)(\{\nu_\alpha(i),\nu_\beta(i)\}\notin E_{\ell_i}).\]
Now, going back to $\bV$, we easily find a condition $q^*\in\bP_{\zeta^*}$
such that $q^*,\alpha,\beta$ are as required.
\end{proof}

Let $q^*,\alpha,\beta$ be as guaranteed by \ref{clx1}. For $j<N$ such that
$\xi^\alpha_j\in I_\ell$, $\ell\leq k$, let $H_j:\dom(\bar{c}(q_\alpha,
\xi^\alpha_j))\longrightarrow\dom(\bar{c}(q_\beta,\xi^\beta_j))$ be the
order preserving bijection (see clause (5) above). We define a condition
$q^+\in\bP_{\gamma}$ as follows:\\
$\dom(q^+)=\dom(q^*)\cup\dom(q_\alpha)\cup\dom(q_\beta)$ and
\begin{itemize}
\item $q^+\rest \zeta^*=q^*$,
\item if $\xi^\alpha_j\in I_\ell$, $n^*\leq j<N$, $\ell\leq k$ then 
\[q^+(\xi^\alpha_j)=\bar{c}_0(q_\alpha,\xi^\alpha_j)\quad\mbox{ and }\quad 
q^+(\xi^\beta_j)=\bar{c}_0(q_\beta,\xi^\beta_j),\]
\item if $n^*\leq j<N$, $\xi^\alpha_j\notin\bigcup\limits_{\ell\leq k}
I_\ell$ then $q^+(\xi^\alpha_j)=q_\alpha(\xi^\alpha_j)$, $q^+(\xi^\beta_j)=
q_\beta(\xi^\beta_j)$. 
\end{itemize}
It should be clear that $q^+\in\bP_\gamma$ is a condition stronger than both
$q_\alpha$ and $q_\beta$. If $i<m$ and $\varepsilon<\omega_1$ then 
\[\nu_\alpha(i)=(\xi^\alpha_j,\varepsilon)\ \Leftrightarrow\ \nu_\beta(i)=
(\xi^\beta_j,H_j(\varepsilon)).\]
If $i\in m\setminus w^*$, $\nu_\alpha(i)=(\xi^\alpha_j,\varepsilon)$,
$n^*\leq j<N$ and $\dot{p}^*_{\varepsilon,\xi^\alpha_j}$,
$\dot{p}^*_{H(\varepsilon),\xi^\beta_j}$ are the names for $\varepsilon^{\rm
th}$ ($H(\varepsilon)^{\rm th}$ respectively) generic reals added by
$\dot{\bQ}_{\xi^\alpha_j}$ ($\dot{\bQ}_{\xi^\beta_j}$, resp.) then 
\[\begin{array}{ll}
q^+\forces_{\bP_\gamma}&\mbox{`` }\dot{p}^*_{\varepsilon,\xi^\alpha_j}\in
\bar{U}_0(q_\alpha,\xi^\alpha_j)(\varepsilon)=\bar{U}_0(q_\beta,\xi^\beta_j)
(H(\varepsilon))\quad\mbox{ and}\\
\ &\ \ \dot{p}^*_{H(\varepsilon),\xi^\beta_j}\in\bar{U}_1(q_\beta,
\xi^\beta_j)(H(\varepsilon))=\bar{U}_1(q_\alpha,\xi^\alpha_j)(\varepsilon)
\mbox{ ''.}
  \end{array}\]
If $i\in w^*$ then look at the choice of $q^*,\alpha,\beta$ (see \ref{clx1}).
Putting everything together we conclude that
\[q^+\forces\mbox{`` }(\forall i<m)(\{\dot{\nu}_\alpha(i),\dot{\nu}_\beta(i)
\}\notin E_{\ell_i})\mbox{ ''},\]
a contradiction
\end{proof}

\section{Proof of Theorem \ref{main}}
Let $\kappa$ be regular cardinal such that $\kappa=\kappa^{<\kappa}\geq
\aleph_2$. By induction on $\xi\leq\kappa$ we build a finite support
iteration $\langle\bP_\xi,\dot{\bQ}_\xi:\xi<\kappa\rangle$ and sequences
$\langle\dot{\bP}^\xi,I_\xi:\xi<\kappa\rangle$, $\langle\dot{A}^\xi_\zeta:
\xi<\zeta\leq\kappa\rangle$ and $\langle\dot{p}^\zeta:\zeta\in I_\xi\rangle$
such that for each $\xi,\xi_0,\xi_1<\kappa$  
\begin{enumerate}
\item[(i)]   $\dot{\bQ}_\xi$ is a $\bP_\xi$--name for a ccc forcing notion
on a bounded subset of $\kappa$, 
\item[(ii)]  $I_\xi\in [\{2\cdot\alpha:\xi<\alpha<\kappa\}]^{\textstyle
\kappa}$, $\dot{\bP}^\xi$ is a $\bP_\xi$--name for an element of $\cK$, and
if $\xi_0\neq\xi_1$ then $I_{\xi_0}\cap I_{\xi_1}=\emptyset$,
\item[(iii)] for $\zeta\in I_\xi$, $\dot{p}^\zeta$ is a $\bP_\zeta$--name for
a condition in $\dot{\bP}^\xi$ for which $\bC(\dot{\bP}^\xi,\dot{p}^\zeta)$
is defined, 
\item[(iv)]  if $\zeta\in I_\xi$ then $\dot{\bQ}_\zeta$ is (equivalent to)
$\bC^{\omega_1}(\dot{\bP}^\xi,\dot{p}^\zeta)$,
\item[(v)]   if $\xi<\zeta\leq\kappa$ then $\dot{A}^\xi_\zeta$ is the
$\bP_\zeta$--name for the set $\bigcup\{\dot{Z}^{\dot{\bP}^\xi}_{
\dot{p}^\varepsilon}:\varepsilon\in I_\xi\cap\zeta\}\subseteq\dot{\bP}^\xi
\times\omega_1$ (where $\dot{Z}^{\dot{\bP}^\xi}_{\dot{p}^\varepsilon}$ is
the generic object added by $\dot{\bQ}_\zeta$; compare \ref{maintech}),
\item[(vi)]  for each $\zeta\leq\kappa$ 
\[\forces_{\bP_\zeta}\mbox{`` }\{(\dot{A}^\varepsilon_\zeta,\cU(
\dot{\bP}^\varepsilon),E(\dot{\bP}^\varepsilon)\cap [
\dot{A}^\varepsilon_\zeta]^{\textstyle 2}): \varepsilon<\zeta\}\mbox{ is 
an $\cS$--family of good graphs '',}\]
\item[(vii)] if $\dot{\bQ}$ is a $\bP_\kappa$--name for a ccc forcing notion
on a bounded subset of $\kappa$ then
\[|\{\zeta<\kappa:\ \ \forces_{\bP_\zeta}\mbox{`` }\dot{\bQ}_\zeta=\dot{\bQ}
\mbox{ ''}\}|=\kappa,\]
\item[(viii)] if $\dot{\bP}$ is a $\bP_\kappa$--name for an element of $\cK$
then then for some $\varepsilon<\kappa$ we have
\[\forces_{\bP_\varepsilon}\mbox{`` }\dot{\bP}^\varepsilon=\dot{\bP}\mbox{
''\ \  and\ \ }\forces_{\bP_\kappa}\mbox{`` the set }\{\dot{p}^\zeta:\zeta
\in I_\varepsilon\} \mbox{ is dense in $\dot{\bP}$ ''.}\]
\end{enumerate}
We use the standard bookkeeping arguments to choose the lists $\langle
\dot{\bP}^\xi,I_\xi:\xi<\kappa\rangle$, $\langle\dot{p}^\zeta:\zeta\in I_\xi
\rangle$ so that clauses (ii), (iii) and (viii) are satisfied. Similarly we
choose a list $\langle\dot{\bQ}'_\xi:\xi<\kappa\rangle$ of all
$\bP_\kappa$--names for partial orders on bounded subsets of $\kappa$ so
that each name appears $\kappa$ many times in the list, and $\dot{\bQ}'_\xi$
is a $\bP_\xi$--name (for $\xi<\kappa$) (this list will be used to take care
of clauses (i), (vii)).  

Now we have to be more specific. So suppose that for some $\xi<\kappa$ we
have already defined the iteration $\langle \bP_\zeta,\dot{\bQ}_\zeta:\zeta<
\xi\rangle$. If $\xi$ is a limit ordinal, before we go further we should
argue that the clause (vi) is satisfied by the limit $\bP_\xi$, i.e., 
\[\forces_{\bP_\xi}\mbox{`` }\cF_\xi\stackrel{\rm def}{=}\{(\dot{A
}^\varepsilon_\zeta,\cU(\dot{\bP}^\varepsilon),E(\dot{\bP}^\varepsilon)\cap
[\dot{A}^\varepsilon_\zeta]^{\textstyle 2}): \varepsilon<\zeta\}\mbox{ is an
$\cS$--family of good graphs ''.}\]
But this is immediate by \ref{maintech} --- if a problem occurs than it is
caused by a finite subfamily of $\cF$ and we may assume that the respective
forcing notions $\bP^{\varepsilon_\ell}$ are from the ground model.

Suppose $\xi=2\cdot\alpha+1$. Then we look at the $\bP_\alpha$--name
$\bQ'_\alpha$ and we ask if, in $\bV^{\bP_\xi}$, it is a ccc forcing
notion. If not that we let $\bQ_\xi$ be the Cohen real forcing. If yes,
then we we ask if (in $\bV^{\bP_\xi}$) it forces that $\cF_\xi$ remains an
$\cS$--family. If again yes, then we let $\dot{\bQ}_\xi$ be
$\dot{\bQ}'_\alpha$; otherwise $\dot{\bQ}_\xi=\bP^{\cF_\xi}$ (see
\ref{notccc}). In any case we are sure that the relevant instances of
clauses (i)--(viii) are satisfied (remember \ref{conclusion}).

Assume now that $\xi=2\cdot\alpha\in I_\zeta$, $\zeta<\kappa$. Then clause
(iv) determines $\dot{\bQ}_\xi$. We should show that the clause (vi) holds
true. Suppose that we may find a condition $q\in\bP_{\xi+1}$, an integer
$m<\omega$, $\varepsilon_0,\ldots,\varepsilon_{m-1}\leq\xi$ and 
$\bP_{\xi+1}$--names $\dot{\nu}_\alpha$ (for $\alpha<\omega_1$) of sequences
of length $m$ such that 
\[\begin{array}{ll}
q\forces_{\bP_{\xi+1}}&\mbox{`` }(\forall\alpha<\beta<\omega_1)(
\dot{\nu}_\alpha\neq\dot{\nu}_\beta)\ \&\ (\forall\alpha<\omega_1)(\forall i<
m)(\dot{\nu}_\alpha(i)\in \dot{A}^{\varepsilon_i}_{\xi+1})\\
\ &\ \ (\forall\alpha<\beta<\omega_1)(\exists i<m)(\{\dot{\nu}_\alpha(i),
\dot{\nu}_\beta(i)\}\in E(\dot{\bP}^{\varepsilon_i}))\mbox{ ''.}
  \end{array}\]
We may additionally demand that for some $k<m$ we have
\[\begin{array}{ll}
q\forces_{\bP_{\xi+1}}&\mbox{`` }(\forall\alpha<\omega_1)(\forall i<k)
(\dot{\nu}_\alpha(i)\in \dot{A}^{\varepsilon_i}_{\xi})\quad\mbox{ and}\\
\ &\ \ (\forall\alpha<\omega_1)(\forall i\in [k,m))(\varepsilon_i=\xi\ \&\
\dot{\nu}_\alpha(i)\in \dot{Z}^{\dot{\bP}^\zeta}_{
\dot{p}^\xi})\mbox{ ''.}
  \end{array}\]
\begin{claim}
\label{clx2}
Suppose that $\bP\in\cK$, $p\in\bP$ (and $\bC(\bP,p)$ is defined). Then 
\[\forces_{\bC^{\omega_1}(\bP,p)}\mbox{`` }(\forall s_0,s_1\in
\dot{Z}^\bP_p)(\{s_0,s_1\}\notin E(\bP))\mbox{ ''.}\]
\end{claim}

\begin{proof}[Proof of the claim] 
Like \ref{oneforce}.
\end{proof}
It follows from \ref{clx2} that 
\[q\forces_{\bP_{\xi+1}}\mbox{`` }(\forall\alpha<\beta<\omega_1)(\exists i<
k)(\{\dot{\nu}_\alpha(i),\dot{\nu}_\beta(i)\}\in
E(\dot{\bP}^{\varepsilon_i}))\mbox{ ''.}\] 
But, by \ref{conclusion}, we have 
\[\forces_{\bP_{\xi+1}}\mbox{`` $\cF_\xi$ is an $\cS$--family of good graphs
'',}\]  
so we get an immediate contradiction. 

Finally if $\xi=2\cdot\alpha\notin\bigcup\limits_{\zeta<\kappa} I_\zeta$
then we let $\bQ_\xi$ be the Cohen real forcing (again all clauses are
preserved).

The construction is complete. We claim that the limit forcing notion
$\bP_\kappa$ is as required in \ref{main}. Clearly it satisfies the ccc and
(a dense subset of it) is of size $\kappa$. Clause (vii) guarantees that
$\forces_{\bP_\kappa}$`` $\con=\kappa\ \&\ \MA$ ''. It follows from the
clause (vi) and \ref{maintech} that 
\[\forces_{\bP_\xi}\mbox{`` }\cF_\kappa\stackrel{\rm def}{=}\{(\dot{A
}^\varepsilon_\zeta,\cU(\dot{\bP}^\varepsilon),E(\dot{\bP}^\varepsilon)\cap
[\dot{A}^\varepsilon_\zeta]^{\textstyle 2}): \varepsilon<\kappa\}\mbox{ is an
$\cS$--family of good graphs ''.}\]
Clauses (v)+(viii) (and the definition of $E(\dot{\bP}^\varepsilon)$) imply
that for each $\varepsilon<\kappa$
\[\forces_{\bP_\kappa}\mbox{`` }(\dot{A}^\varepsilon_\zeta,\cU(\dot{\bP
}^\varepsilon),E(\dot{\bP}^\varepsilon)\cap [\dot{A}^\varepsilon_\zeta]^{
\textstyle 2})\mbox{ is densely representably by
$\bot_{\dot{\bP}^\varepsilon}$ ''}.\]
Consequently, by \ref{notFAP} and clause (viii) we get 
\[\forces_{\bP_\kappa}\mbox{`` }(\forall\bP\in\cK)(\neg\FA_{\aleph_1}(\bP))
\mbox{ '',}\]
finishing the proof of \ref{main}.

\begin{corollary}
It is consistent with $\MA+\neg{\rm CH}$ that any forcing notion $\bP\in\cK$
collapses $\con$ to $\omega_1$ (and thus is not $\omega_1$--proper).
\end{corollary}

\section{Open problems}
The model constructed in the previous section provides
$(\forall\bP\in\cK)(\neg\FA_{\aleph_1}(\bP))$ by dealing with each forcing
$\bP\in\cK$ separately. We would like to have one common reason for
$\neg\FA_{\aleph_1}(\bP)$ for all forcing notions $\bP\in\cK$, i.e., a
combinatorial principle $\cP$ which is consistent with $\MA+\neg {\rm CH}$
and which implies $(\forall\bP\in\cK)(\neg\FA_{\aleph_1}(\bP))$. A possible
candidate for a principle like that was already pointed in \cite[\S
2]{CRSW93}. As we stated in the Introduction, our method is a slight
generalization of that of Steprans. Steprans' method in turn was based on
the proof of Abraham, Rubin and Shelah \cite{ARSh:153} that it is consistent
with $\MA+\neg {\rm CH}$ that there are two non-isomorphic $\aleph_1$--dense
sets of reals. In the latter proof, a 2--entangled set of reals was
used. This leads us to the following question 
\begin{problem}
[Compare {\cite[Question 2.4]{CRSW93}}]
Does the existence of a 2--entangled set of reals of size $\aleph_1$ imply 
$(\forall\bP\in\cK)(\neg\FA_{\aleph_1}(\bP))$?
\end{problem}

If one tries to repeat the proof of \cite[Theorem 2.1]{JMSh:372} for
elements of $\cK$ then one gets into some problems in cases 1,3 of
Definition \ref{classK}. Possible reason for it is that a proof as in
\cite[Theorem 2.1]{JMSh:372} would give a property that seems to be stronger
than $\neg\FA_{\aleph_1}(\bP)$.

\begin{definition}
Let $\bP$ be a forcing notion of size $\con$, $\bar{p}=\langle p_i: i<\con
\rangle\subseteq\bP$. We say that $\bar{p}$ is an JMSh--sequence if 
\begin{enumerate}
\item[$(\oplus)_{\rm JMSh}$] given $\langle F_\alpha:\alpha<\omega_1\rangle$
pairwise disjoint finite subsets of $\con$, there exist $\alpha<\beta<
\omega_1$ such that
\[(\forall i\in F_\alpha)(\forall j\in F_\beta)(p_i\, \bot_{\bP}\; p_j).\]
\end{enumerate}
\end{definition}

\begin{proposition}
\label{JMShseq}
Suppose that $\bP$ is a forcing notion of size $\con$ such that
\begin{enumerate}
\item above any condition in $\bP$, there is an antichain (in $\bP$) of size
$\con$, and
\item there is an JMSh--sequence $\bar{p}\subseteq\bP$ which is dense in
$\bP$.
\end{enumerate}
Then $\neg\FA_{\aleph_1}(\bP)$.
\end{proposition}

\begin{problem}
\begin{enumerate}
\item Is the existence of dense JMSh--sequences in $\bP$ really stronger
than $\neg\FA_{\aleph_1}(\bP)$ (for $\bP$ of size $\con$ satisfying the
assumption \ref{JMShseq}(1)) ?
\item Is it consistent with $\MA+\neg{\rm CH}$ that for each $\bP\in\cK$
there is a dense JMSh--sequence in $\bP$?
\end{enumerate}
\end{problem}

On the other hand, Judah, Miller and Shelah \cite{JMSh:372} and Goldstern,
Johnson and Spinas \cite{GJS94} showed that ${\bf MA}_{\omega_1}(\mbox{ccc}
)$ implies the forcing axiom for the Miller and Laver forcing notions. This
gives a strong expectation that ${\bf MA}_{\omega_1}(\mbox{ccc})$ implies
forcing axioms for most of forcing notions (with norms) adding unbounded
reals. Brendle \cite[Proposition 5.1]{Br96a} showed that $\MA$ implies that
the Laver forcing, the Mathias forcing, the Miller forcing and the
Blass-Shelah forcing are $\alpha$--proper for all $\alpha<\con$. Again, one
expects that this could be generalized further.

\begin{problem}
Let $\cK^\bot$ be the class of the forcing notions of \cite{RoSh:470} which
are not in $\cK$ for nontrivial reasons.
\begin{enumerate}
\item Does $\MA+\neg{\rm CH}$ imply $\FA_{\aleph_1}(\bP)$ for all $\bP\in
\cK^\bot$ ?
\item Does $\MA+\neg{\rm CH}$ imply that all $\bP\in\cK^\bot$ are
$\alpha$--proper (for all $\alpha<\con$) ?
\end{enumerate}
\end{problem}

Finally, possible further generalizations of the present paper could go into
the direction of nep/snep forcing notions of Shelah \cite{Sh:630},
\cite{Sh:669}. 

\begin{problem}
\begin{enumerate}
\item Is $\MA+\neg{\rm CH}$ consistent with the failure of
$\FA_{\aleph_1}(\bP)$ for all snep $\baire$--bounding forcing notions $\bP$
which do not have ccc above any condition?
\item Does $\MA+\neg{\rm CH}$ imply $\FA_{\aleph_1}(\bP)$ for all snep
forcing notions $\bP$ adding unbounded reals?
\end{enumerate}
\end{problem}


\end{document}